\documentclass[12pt]{article}
\begin{document}
\begin{center}
\large{\bf{Random symmetrizations of measurable sets}}
\bigskip

{\large{\bf{ Aljo\v{s}a Vol\v{c}i\v{c}}}}

\end{center}
\bigskip

\begin{abstract}
In this paper we prove almost sure convergence to the ball, in the Nikodym metric, of sequences of random Steiner symmetrizations of bounded Caccioppoli and bounded measurable sets, paralleling a result due to Mani-Levitska concerning convex bodies.
\smallskip

\end{abstract}

\section {\bf Introduction} 

In [M] Mani-Levitska studied sequences of random Steiner symmetrizations of convex bodies $K\subset I\!\!R^d$, proving that almost surely they converge in the Hausdorff distance to the ball centered at the origin and having the same volume. His result improved an old theorem due to Gross [G] who proved that, given a convex body $K$, there exists a sequence of directions such that its successive Steiner symmetrizations converge to that ball.  Since Gross' theorem holds more generally for compact sets, Mani-Levitska conjectured that his extension can also be proved for the class of compact sets.

This paper tackles and solves the analogous question within the class of bounded measurable sets,  where the natural metric is the symmetric-difference distance (called also Nikodym distance). The results seem to be of independent interest, but we believe that they will bring to the solution of Mani-Levitska's question, which is one of the motivations of this paper.

In Section 2 we will prove some preliminary results concerning the Steiner symmetrization and the central moment of inertia.

Section 3 contains the main results,  showing first that the desired conclusion holds for bounded
Caccioppoli sets. The extension to the class of bounded measurable sets is then obtained by approximation. 

The last section is devoted to numberless open questions suggested by the present results and some problems collected from the pertaining literature.

\section {\bf Notations, definitions and preliminary results}

We denote by $\lambda_d$ the $d$-dimensional Lebesgue measure and by  $\cal M$ the family of all bounded measurable sets of the $d$-dimensional euclidean space $I\!\!R^d$, which will be called {\it summable} sets. Two summable sets $A$ and $B$ are {\it equivalent} ($A\sim B$) if the Lebesgue measure of their symmetric difference  $A\bigtriangleup B$ is zero, otherwise we say that they are essentially different. Whenever necessary, we will adopt the usual and useful ambiguity between summable sets and the corresponding equivalence classes.

We shall say that a subset $\cal A \subset \cal M$ is {\it bounded} if there is a bounded measurable set $G$ containing all members of $\cal A$.

On $\cal M$ we define the following pseudo metric, called the Nikodym distance or also the symmetric-difference distance.   If $A$ and $B$ are summable sets, their pseudo-distance is 
$$d_N(A,B)=\lambda_d (A\bigtriangleup B)\,.$$

The quotient space $\cal M_{/\sim}$ is a complete metric space with respect to $d_N$.
\smallskip

By $\lambda^*$ we shall denote the outer one-dimensional measure: if $E\subset I\!\!R$ is any set, $\lambda^*(E)$ is the smallest measure of a measurable set containing $E$.
\smallskip

A set $R=\prod_{i=1}^d [a_i, b_i]$ will be called a rectangle. The origin of $I\!\!R^d$ will be denoted by $o$.
By $B(q,\rho)$ we denote the closed ball centered in $q$ having radius $\rho$. The unit sphere, i.e. the set of all unit vectors, will be denoted by $S^{d-1}$. The volume of the unit ball of $I\!\!R^d$ centered at the origin, $B(o,1)$,  is traditionally denoted by $\kappa_d$. 

If $A$ is a summable set, then a ball having the same volume has radius $\rho(A)=\left(\frac{ \lambda_d(A)}{\kappa_d}\right)^{\frac{1}{d}}$.
\smallskip

The scalar product between two vectors $u$ and $v$ of $I\!\!R^d$ will be denoted by $u\cdot v$.
\bigskip

We shall now define a fundamental notion, introduced (for convex sets) by Jacob Steiner [S] in order to solve the classical isoperimetric problem. His proof contained  an inevitable gap (for that time, since the concept of compactness was not available) which has been fixed only much later, but his beautiful geometric idea is still an indispensable tool in many branches of geometry and analysis.
\bigskip

\noindent {\bf Definition 2.1} Let $A$ be a measurable set in $I\!\!R^d$, $u$ a unit vector and let $l_u$ be the line through the origin parallel to $u$. We denote by $u^{\perp}$ the hyperplane orthogonal to $u$. For each $x\in u^{\perp}$, let $c(x)$ be defined as follows. If $A\cap (l_u +x)$ is empty, let $c(x)=\emptyset$. Otherwise, let $c(x)$ be the possibly degenerate (and possibly infinite) closed segment on $l_u+x$ centered at $x\in u^{\perp}$ whose length is equal to $\lambda^*(A\cap (l_u +x))$. 

The union of all the line segments $c(x)$ is called the {\it Steiner symmetral} of $A$ and will be denoted by $S_uA$. The mapping $S_u$ from the family of measurable sets into itself is called {\it Steiner symmetrization}.
\bigskip

In the literature one can find several definitions of Steiner symmetral of a measurable set. The one we use can be found in [AFP].  It differs slightly from the one in  ([Ga], Definition 2.1.3), where it is assumed that $c(x)=\emptyset$ when $A\cap (l_u +x)$ is not measurable. With our definition, Steiner symmetrization is monotone in the sense that if $A\subset B$, then $S_uA\subset S_uB$.

The symmetral $S_uA$ is measurable, summable, Borel, compact or convex when $A$ is measurable, summable, Borel, compact or convex (compare [Ga], Theorem 2.1.4) .
\smallskip

Observe also that Steiner symmetrization is compatible with the equivalence relation in $\cal M$.
\bigskip

The next lemmas are presented in preparation of Section 3 and mostly they are not stated in full generality, since we do not need to look at classes of sets beyond $\cal M$.
\bigskip

\noindent {\bf Lemma 2.2} {\it The Steiner symmetrization $S_u$ is Lipschitz with constant $1$ on $\cal M_{/\sim}$. }
\bigskip
 
\noindent {\it Proof.} By the Fubini theorem
$$\lambda_d (S_uA \bigtriangleup S_uB)=\int_{u^{\perp} }\left(\int_{l_u+x} \chi_{(S_uA \bigtriangleup S_uB)}(t)\,dt\right)\,d\lambda_{d-1}(x)$$
$$\le \int_{u^{\perp} }\left(\int_{l_u+x} \chi_{(A \bigtriangleup B)}(t)\,dt\right)\,d\lambda_{d-1}(x)=\lambda_d (A \bigtriangleup B)\,.$$

\smallskip
\noindent since almost everywhere
$$\lambda_1(S_uA\cap (l_u+x) \bigtriangleup S_uB\cap (l_u+x) )\le \lambda_1(A\cap (l_u+x)  \bigtriangleup B\cap (l_u+x) )\,.\eqno {\bf\star}$$

Directions in $I\!\!R^{d}$ will be identified by unit vectors, hence by elements of the unit sphere $S^{d-1}$. Since $u$ and $-u$ define the same direction, we shall consider on $S^{d-1}$ the following equivalence relation: $u\sim v$ if and only if $u=v$ or $u=-v$. On $S^{d-1}/{\sim}$ we consider the distance defined by
$$d(u,v)= \min (\|u-v\|, \|u+v\|)\,.$$

The probability on the set of directions  $S^{d-1}/{\sim}$ is defined by the normalized surface area restricted to the $\sigma$-algebra of symmetric Borel sets, i.e. to those Borel sets $D\subset S^{d-1}$ such that if $u\in D$, then $-u\in D$.
\bigskip

\noindent {\bf Lemma 2.3} {\it If $A\subset I\!\!R^d$ is a summable set, then the mapping 
$\varphi  : u \mapsto S_u A$
 is uniformly  continuous  on $S^{d-1}/\sim$.}
\bigskip

\noindent {\it Proof.} Since $S^{d-1}/\sim$ is compact, we only need to prove continuity. 

The conclusion is obvious if $A$ is a rectangle or a finite union of rectangles. Otherwise, given $\varepsilon>0$, let $R_{\varepsilon}$ be a finite union of rectangles such that $d_N(A, R_{\varepsilon})<\frac{\varepsilon}{3}$. There exists a $\delta>0$ such that if $d(u, u_0)<\delta$, then $d_N (S_uR_{\varepsilon}, S_{u_0}R_{\varepsilon})<\frac{\varepsilon}{3}$.

We have therefore that 
$$d_N(S_u A, S_{u_0}A)\le d_N(S_u A, S_{u}R_{\varepsilon})
+d_N(S_u R_{\varepsilon}, S_{u_0}R_{\varepsilon})+d_N(S_{u_0} R_{\varepsilon}, S_{u_0}A)$$
and the conclusion follows since, by Lemma 2.2, $d_N(A, R_{\varepsilon})<\frac{\varepsilon}{3}$ implies that  $d_N(S_v A, S_v R_{\varepsilon})<\frac{\varepsilon}{3}$ for any $v\in S^{d-1}/\sim$.$\,\,\,{\bf\star}$
\bigskip

A random sequence of symetrizations corresponds to a random sequence of directions  choosen independently and uniformly in $S^{d-1}/\sim$, i.e. to an element $U=(u_1, u_2, \dots, u_n, \dots )$ of the probability space
$$\Omega= \prod_{n=1}^{\infty}(S^{d-1}/{\sim})_n\,,$$ 
where $(S^{d-1}/\sim)_n =S^{d-1}/{\sim}$ for every $n$ and the probability is the ordinary product probability.
\smallskip

By $U_n$ we denote the $n$-tuple $(u_1, u_2, \dots, u_n )$ and by $S_{U_n}A$ the set obtained symmetrizing successively $A$ with respect to $u_1, u_2, \dots, u_n $, i.e. $S_{U_n}A=S_{u_n}S_{u_{n-1}}\dots  S_{u_1}A$. Note that the order is important, since Steiner symmetrizations are, in general, not commutative.
\bigskip

The following proposition is the analogue, in this setting, of Lemma 1 of [M].
\bigskip

\noindent {\bf Proposition 2.4} {\it Let $A$ be a summable set and $\{ u_n\}$ be a random sequence of independent  directions.  If we denote by $b_n$ the barycenters of $ S_{U_n}A $, the successive Steiner symmetrizations of $A$, then almost surely the sequence $\{ b_n\}$ converges to the origin.}
\bigskip

\noindent {\it Proof.} Since $b_{n+1}$ is the projection of $b_n$ on the hyperplane $u_{n+1}^{\perp}$, we have $\|b_n\|^2 =\|b_n-b_{n+1}\|^2+\|b_{n+1}\|^2\,,$
and therefore, for any $n$, $\|b_{n+1}\|\le \|b_n\|$. If $b_n=o$ for some $n$, there is nothing to prove. Let otherwise $w_n=\frac{b_n}{\|b_n \|}$. If $|w_n \cdot u_{n+1}|\ge \frac{1}{\sqrt{2}}$, then $\|b_{n+1}\|\le \frac{1}{\sqrt{2}}\|b_n\|$. 

Denote by $\zeta_d$ the probability (which depends on the dimension $d$ but not on $w$) of the double cap $\{u: | u\cdot w |\ge \frac{1}{\sqrt{2}}\}$.
This is the probability that $|w_n \cdot u_{n+1}|\ge \frac{1}{\sqrt{2}}$. Hence, by the (divergence part of the) Borel-Cantelli lemma, almost surely $\|b_{n+1}\|\le \frac{1}{\sqrt{2}}\|b_n\|$ for infinitely many indices $n$, and the conclusion follows. ${\bf\star}$
\bigskip

\noindent {\bf Definition 2.5} Given a summable set $A$, its {\it (central) moment of inertia} is defined by
$$\mu(A)=\int_A \|z\|^2\, d\lambda_d(z)\,.$$
\smallskip

Observe that $\mu$ is defined for all summable sets and that if $A$ and $B$ are equivalent, then $\mu(A)=\mu(B)$. Put $\mu_d=\mu(B(o,1))$.

We will now see some useful properties of this functional. 
\bigskip

\noindent {\bf Lemma 2.6} {\it The moment of inertia $\mu$ is uniformly continuous on bounded subsets of $\cal M$. }
\bigskip

\noindent {\it Proof.} Let $\cal A$ be a bounded subset of $\cal M$ and let $G$ be a summable set containing all members of $\cal A$. Since $G$ has finite Lebesgue measure and $\mu$ is absolutely continuous with respect to $\lambda_d$, for every $\varepsilon >0$ there exists a $\delta>0$ such that if $E\subset G$ and $\lambda_d(E)<\delta$, then $\mu(E)<\varepsilon$.

Therefore for any $A_1, A_2\in \cal A$, if $d_N(A_1,A_2)<\delta$ then
$$|\mu(A_1)-\mu(A_2)| \le \mu(A_1\bigtriangleup A_2)<\varepsilon\,,$$
and the conclusion follows.$\,\,\,{\bf\star}$
\bigskip

\noindent {\bf Lemma 2.7} {\it The functional $(u,A)\mapsto \mu (S_u A)$ is uniformly continuous on $S^{d-1}/\sim \times \, \cal A\,$, whenever $\cal A$ is a bounded subset of $\cal M$. }
\bigskip

\noindent {\it Proof.} Since 
$$|\mu(S_uA)-\mu( S_vB)| \le| \mu(S_uA)-\mu( S_uB)|+|\mu(S_uB)-\mu( S_vB)|\,,$$
 the conclusion follows from Lemma 2.6, using Lemmas 2.2 and 2.3.$\,\,\,{\bf\star}$
\bigskip

\noindent {\bf Lemma 2.8} {\it If $A$ is any measurable set which is essentially different from $B=B(o,1)$ and $\lambda_d(A)=\kappa_d$, then $ \mu(A)> \mu_d$. }
\bigskip

\noindent {\it Proof.} We have
$$\mu_d=\mu (B)=\int_B \|z\|^2 \, d\lambda_d =\int_{B\cap A}\|z\|^2  \,d\lambda_d  +  \int_{B\setminus A}\|z\|^2  \,d\lambda_d  $$
$$<\int_{B\cap A}\|z\|^2  \,d\lambda_d +  \int_{A\setminus B}\|z\|^2  \,d\lambda_d  =\mu(A)\,.$$
The strict inequality is due to the fact that $B\setminus A$ and $A\setminus B$ have the same positive measure and  $\|z\|^2 <1$ almost everywhere on $B\setminus A$, while  $\|z\|^2 >1$ almost everywhere on $A\setminus B$.$\,\,\,{\bf\star}$
\bigskip

Note that the integral representing $\mu(A)$ may diverge. The conclusion can also be refrased saying that the unit ball is the unique minimizer of $\mu$ among all the measurable sets having measure $\kappa_d$.
\bigskip

The next lemma will show that Steiner symmetrization in appropriate directions strictly diminishes the moment of inertia of any set in $\cal M$ which is not equivalent to a ball centered at the origin.
\bigskip

\noindent {\bf Lemma 2.9} {\it If $A\in \cal M$ is essentially different from the ball $B=B(o,\rho(A))$, then there exist a direction $v$ and a positive $\delta$ such that
$$\mu(S_u A)<\mu(A)$$
for all $u$ such that $d(u,v)<\delta$. }
\bigskip

\noindent {\it Proof.} By assumption, the sets $E=A\setminus B$ and $F=B\setminus A$ have the same positive measure. Therefore there exist, by the Lebesgue density theorem, points $q_2$ and $q_1$ which are of density $1$ for $E$ and $F$, respectively. Then for every $\varepsilon >0$ there exists a $\rho>0$ such that, if we put $B_2=B(q_2,\rho)$ and $B_1=B(q_1,\rho)$, we have  $B_2\cap B= \emptyset$ and $B_1 \subset B$ and moreover
$$\lambda_d (A\cap B_2)>(1-\varepsilon) \lambda_d (B_2) \mbox{\,\,\,\,\,and\,\,\,\,\,} \lambda_d ( B_1\cap A)<\varepsilon \lambda_d (B_1)\,.$$

We choose a coordinate system on $I\!\!R^d$ such that $v=(0,0, \dots, 0, 1)$ and denote a generic point of $I\!\!R^d$ as  $z=(x,t)$, with $x\in I\!\!R^{d-1}$ and $t\in I\!\!R$. 

Let $q_2=(0,0, \dots, t_2)$ and $q_1=(0,0, \dots, t_1)$. Since $q_2$ does not belong to $B$, while $q_1$ is interior to $B$, we have that $t_2>|t_1|$.

We will assume $\varepsilon <\frac{1}{8}$ and we are of course allowed to take $0<\rho<\frac{t_1+t_2}{8}$.

Let $v=\frac{q_2-q_1}{\|q_2-q_1\|}$ and symmetrize with respect to $v$. Then we have
$$\mu (A)=\int_{v^{\perp}}\left( \int _{A\cap (l_v+x)} \|z\|^2 \,dt \right)\, d\lambda_{d-1}(x)$$
$$=\int_{v^{\perp}}\left( \int _{A\cap (l_v+x)} t^2 \,dt + \int _{c(x)} \|x\|^2 \,dt \right)\, d\lambda_{d-1}(x)\,.$$

\smallskip 
\noindent The last equality holds because $\|x\|^2$ does not depend on $t$. It follows that any variation of $\mu(S_v A)$ with respect to $\mu(A)$ depends only on the first integral in the brackets.

Consider now the set $A'=(A\setminus(B_1\cup B_2))\cup ((B_2\cap A-(t_2-t_1)v))\cup ((B_1\cap A+(t_2-t_1)v))$ which is obtained from A changing the positions of $B_1\cap A$ and $B_2 \cap A$. Clearly $S_v A=S_v A'$ and hence $\mu(S_v A)=\mu(S_v A')$. 

Let us show that $\mu(A')<\mu(A)$.

Since $A$ and $A'$ differ only in $B_1\cup B_2$, we have
$$\mu(A)-\mu(A')=\int_{B_2\cap A}t^2\,d\lambda_d +\int_{B_1\cap A}t^2\,d\lambda_d $$
$$-\left(\int_{B_2\cap A}(t-s)^2\,d\lambda_d +\int_{B_1\cap A}(t+s)^2\,d\lambda_d \right)\,,\eqno (1)$$
where $s=t_2-t_1$. Expanding the squares in (1) and simplifying, we may rewrite the previous expression as
$$\int_{B_2\cap A}(2ts-s^2)\,d\lambda_d -\int_{B_1\cap A}(2ts+s^2)\,d\lambda_d$$
$$=2s\left( \int_{B_2\cap A}(t-\frac{s}{2})\,d\lambda_d -\int_{B_1\cap A}(t+\frac{s}{2})\,d\lambda_d\right)\,.$$

\smallskip
Observe now that $t-\frac{s}{2}\ge t_2-\rho -\frac{s}{2}$ on $B_2 \cap A$ and $t+\frac{s}{2}\le t_1+\rho + \frac{s}{2}$ on $B_1 \cap A$, so if we use the bounds we assumed for $\varepsilon$ and $\rho$, the last expression is  bounded from below by 
$$2s\left( \lambda_d (B_2\cap A)( t_2-\rho -\frac{s}{2} ) -\lambda_d (B_1\cap A)( t_1+\rho + \frac{s}{2})\right)$$
$$\ge 2s\left( \lambda_d (B_2)( t_2-\rho -\frac{s}{2})(1-\varepsilon) -\lambda_d (B_1)( t_1+\rho + \frac{s}{2})\varepsilon \right)$$
$$= 2s \lambda_d (B_2)\left( \frac{t_1+t_2}{2}- \varepsilon (t_1+t_2)  -\rho \right) >\frac{s}{2}\lambda_d (B_2)(t_1+t_2)>0$$
It follows therefore that
$$\mu(S_v A)=\mu(S_v A')\le\mu(A')<\mu(A)\,. \eqno (2)$$

To conclude the proof, let us recall that, by Lemma 2.7, the functional $u\mapsto \mu(S_uA)$ is uniformly continuous on $S^{d-1}/\sim$, therefore there exists a $\delta >0$ such that the inequality (2) is preserved for all $u$ such that $d(u,v)<\delta$.$\,\,\bf \star$

\section {\bf Main results}  

This section is devoted to random Steiner symmetrizations of Cacciopoli and summable sets.
\bigskip

\noindent {\bf Definition 3.1} A measurable subset $C$ of $I\!\!R^d$ is called a {\it Caccioppoli set} if 
$$\lambda_d ((C+h)\bigtriangleup C)\le p \cdot |h|\,$$
for some costant $p$ and for every $h\in I\!\!R^d$.
\bigskip

This definition is due to Caccioppoli [C] and has been extensively studied and used by De Giorgi ([DG] and subsequent papers). We will denote by $\cal C$ the family of all bounded Caccioppoli sets. 

For $C\in \cal C$ we will denote by $p(C)$ the De Giorgi-Caccioppoli perimeter of $C$ ([T], p. 85) which can be defined for instance by
$$p(C)=\inf \{ \liminf p(E_n)\}\,,$$
where $\{E_n\}$ is a sequence of smooth subsets of $I\!\!R^d$ such that $d_N(C, E_n)\rightarrow 0$, the infimum is taken over all such sequences, and the perimeter for smooth sets is understood in the ordinary sense.
\bigskip

We will need the following properties on Caccioppoli sets and the perimeter, which can all be found in [T]:
\bigskip

(i) The Steiner symmetral of a Caccioppoli set is a Caccioppoli set. 
\smallskip

(ii) The Steiner symmetrization does not increase the perimeter of a Caccioppoli set. 
\smallskip

(iii) If  $\cal F \subset \cal C$  is a collection of Caccioppoli sets which are contained in a bounded set and have uniformly bounded perimeters, then $\cal F$ is relatively compact with respect to the Nikodym distance. 
\smallskip

(iv) The perimeter is lower semicontinuous on $\cal C$.
\bigskip

Now we prove the key lemma on which the main result is based.
\bigskip

\noindent {\bf Lemma 3.2} {\it Fix $\rho_0>1$, $p_0>0$ and $\varepsilon_0>0$ and consider the family 
${\cal F}={\cal F}(\rho_0, p_0, \varepsilon_0)$ of all Caccioppoli sets such that $\lambda_d(F)=\kappa_d$, contained in $B(o,\rho_0)$ and having perimeters bounded by $p_0$.  Suppose moreover that 
$$\mu(F)\ge \mu_d + \varepsilon_0 \eqno (3)$$
for every $F\in \cal F$.
Then there exist a $\delta_0>0$  and for each $F\in {\cal F}$ a direction $v_F$ such that
$$\mu(S_vF)<\mu(F)-\delta_0$$
for every $v$ such that $d(v,v_F)<\delta_0$.}
\bigskip

\noindent {\it Proof.} Observe first that $\cal F$ is closed and hence compact. We shall prove the lemma by contradiction. Suppose the conclusion is not true. Then, for every $n\in I\!\!N$ there exist $F_n\in {\cal F}$ and, for every direction $v$, a corresponding direction $v_n$ such that
$$d(v_n,v)\le \frac{1}{n}$$
and 
$$\mu(S_{v_n}F_n)\ge \mu(F_n)-\frac{1}{n}\,.$$
Compactness implies that there exists a subsequence $\{F_{n_k}\}$ converging to a Caccioppoli set $F\in \cal F$.

By continuity of the Steiner symmetrization and of the moment of inertia,
$$\mu(S_vF)\ge \mu(F)$$
for every $v\in S^{d-1}$, but since the opposite inequality is obvious, we have $\mu(S_vF)= \mu(F)$ for every direction $v$ and hence, by Lemma 2.9, $F$ is equivalent to $B(o,1)$ which therefore has to belong to $\cal F$. But this contradicts (3).$\,\,\,\star$
\bigskip

\noindent {\bf Remark 3.3} A consequence of Lemma 3.2 is that if $F\in {\cal F}$, there is a $\delta_0 >0$, such that with probability $P_0=P(\{u: d(u,u_0)<\delta_0\})$ the moment of inertia of $F$ will be diminished by $\delta_0$ by a random simmetrization. Note that  this probability depends on $\delta_0$ and hence on $\rho_0$, $p_0$ and $\varepsilon_0$, but does not depend on $F$.
\bigskip

We are now ready to prove the main result of this section.
\bigskip

\noindent {\bf Theorem 3.4} {\it If $F$ is any bounded Caccioppoli set, then with probability $1$ its successive random Steiner symmetrizations $F_n=S_{U_n}F$ converge, with respect to the Nikodym distance, to the ball $B(o,\rho(F))$. }
\bigskip

\noindent {\it Proof.} Since the problem is invariant by dilation, we may suppose that $\lambda_d(F)=\kappa_d$.

Note that for any bounded Caccioppoli set $F$, $\{\mu(F_n)\}$ is a decreasing sequence which tends to $\mu_d$ if and only if $\{F_n\}$ tends to $B=B(o,1)$. This is because the sequence $\{F_n\}$  is relatively compact and all we have to check is that any convergent subsequence converges to $B$. Suppose the contrary, and let a subsequence $\{F_{n_k}\}$ converge to $C$, essentially different from $B$. Then $\{\mu(F_{n_k})\}$ tends to $\mu(C)>\mu_d$,  a contradiction by Lemma 2.8.

If we assume that the conclusion of the theorem does not hold, there exist a bounded Caccioppoli set $F_0$ and a set of positive probability ${\cal U}\subset \Omega$ such that the sequence $\{S_{U_n}F_0\}$ does not converge for any $U\in \cal U$. Therefore there exists an $\varepsilon_0>0$ such that
$$\mu(S_{U_n}F_0)\ge \mu_d + \varepsilon_0 $$ 
for all $n\in I\!\!N$ and all $U\in \cal U$.

 Take $\rho_0>1$ such that $F_0\subset B(o,\rho_0)$. All the successive Steiner symmetrizations $F_n$ of $F_0$ have the same measure and are contained in the same ball. Moreover, the perimeters of the sets $F_n$ are bounded by the perimeter $p_0$ of $F_0$. 

By Lemma 3.2 there exists a $\delta_0>0$ and a probability $P_0$ depending on $\varepsilon_0$, $p_0$ and $\rho_0$  such that 
$$\mu(F_{n+1})<\mu(F_n)-\delta_0 \eqno (4)$$
happens with probability (at least) $P_0$ while with positive probability 
$$\mu(F_n)\ge \mu_d+\varepsilon_0\,. \eqno (5)$$
But the divergence part of the Borel-Cantelli lemma assures that if (5) holds, (4) will almost surely happen infinitely often and therefore, with probability 1, $\mu(F_n)<\mu_d+\varepsilon_0$ for $n$ sufficiently large, a contradiction. $\,\,\,\star$
\bigskip

We shall see now that the same conclusion holds more generally for summable sets.
\bigskip

\noindent {\bf Theorem 3.5} {\it If $A$ is a summable set, then with probability $1$ its successive random Steiner symmetrizations $A_n=S_{U_n}A$ converge, with respect to the Nikodym distance, to the ball $B(o,\rho(A))$. }
\bigskip

\noindent {\it Proof.} Since $A$ is summable, it is contained in a ball $B(o,\rho)$ and the same is true for all its subsequent symmetrizations. We may again assume that $\lambda_d(A)=\kappa_d$. 

Fix $\varepsilon >0$ and let $R$ by a finite union of rectangles (and hence a Caccioppoli set) contained in $B(o,\rho)$ such that $d_N(A,R)<\frac{\varepsilon}{2}$ and $\lambda_d(R)=\kappa_d$.  By Lemma 2.2 we have also $d_N(A_n, S_{U_n}R)<\frac{\varepsilon}{2}$ for any $n$, and hence
$$d_N(A_n,B(o,1))\le d_N(A_n, S_{U_n}R)+d_N( S_{U_n}R, B(o,1))$$
$$<\frac{\varepsilon}{2}+d_N( S_{U_n}R, B(o,1))\,,$$
for any $n\in I\!\!N$.

On the other hand by Theorem 3.4, with probability 1, 
$$\lim_{n\rightarrow \infty}d_N(S_{U_n}R, B(o,1))=0\,,$$
and therefore with probability 1, for all $n$ sufficiently large,
$$d_N(S_{U_n}A, B(o,1))<\varepsilon\,,$$
proving so the claim. $\,\,\,\star$
\bigskip

\noindent {\bf Remark 3.5} The conclusion of the previous theorem holds also if we assume that the random directions are possibly dependent, but pairwise independent ([Ch], Theorem 4.2.5), since the convergence part of the Borel-Cantelli lemma holds also under these weaker conditions.

A similar generalization can be done with Proposition 2.4.

\section {\bf Conclusions and open problems} 

A consequence of Theorem 3.4 is that, given a compact set $K$ of positive measure,  random Steiner symmetrizations converge with probability 1 {\it in the Nikodym distance} to the ball $B(o,\rho(K))$. So we got a different answer to Mani-Levitska's question, since the natural metric on $\cal K$, the collection of compact sets, is the Hausdorff distance, which induces a topology which is not comparable on $\cal K$ with the topology induced by the Nikodym (pseudo-) metric.

We believe however that the methods developed in this paper, combined with some more ideas we have in mind, will allow us to solve the question which originally motivated this paper.
\smallskip 

An interesting but technical question is whether the boundedness of the Caccioppoli and measurable sets can be removed. The obstacle lies clearly in the fact that we use as an essential tool the moment of inertia, which is defined on all summable sets but not on all measurable sets of finite measure. We could consider the class of (possibly unbounded) ${\cal L}_2$ Caccioppoli sets, using finer versions of Lemmas 2.6 and 2.7, since we apply them only on families of sets which are obtained by successive Steiner symmetrizations. But is it possible to do better, for instance using the lower semicontinuity of the perimeter instead of the continuity of the moment of inertia?
\smallskip

A question suggested by Michele Gianfelice is whether the uniform distribution on $S^{d-1}/\sim$ is essential or if one could replace it with other symmetric probability distributions. 

We are inclined to believe that in fact any symmetric distribution leads to the same result.
\smallskip

Another question is whether (pairwise) independence is necessary. We conjecture that pairwise negative correlation (a condition we exploited recently in another context) is a reasonable condition to look at.
\smallskip

The problem admits also some interesting deterministic variants: does $\{S_{U_n}A\}$ converge to $B(o,\rho(A))$ when $\{u_n\}$ is a uniformly distributed sequence on $S^{d-1}/\sim$ (see [KN], Chapter 3, for the definition)? It should be noted that, with probability $1$, a random sequence is uniformly distributed (see [KN], Chapter 3, Theorem 2.2), so there is a considerable overlap between the two conditions. However an open question remains (with the obvious variants): does there exist a uniformly distributed sequence of directions $\{u_n\}$ and a bounded Caccioppoli (convex, compact, summable, having finite measure) set $F$ such that $\{S_{U_n}F\}$ does not converge to $B(o, \rho(F))$ in the appropriate metric?

We believe that the answer to the above question is negative and in fact we conjecture that a much stronger result holds, namely that the density alone of the sequence of directions  $\{u_n\}$ it sufficient for the convergence of $\{S_{U_n}A\}$ in the various settings mentioned above.
\smallskip

Another problem is whether there exists a finite number of directions such that alternating the symmetrization with respect to them we have that $\{S_{U_n}A\}$ is convergent to $B(o,\rho(A))$, with all the variants described above concerning the type of sets (convex, compact, Caccioppoli, summable, having finite measure), and the metric. Some encouraging examples suggest that this may be an interesting question to investigate on.
\smallskip

Several authors studied bounds for the number of successive Steiner symmetrizations required to transform any convex body $K$ of given volume to a convex body which is ``close" to the ball having the same volume ([H], [BLM], [BG], [KM1], [KM2]). It would be interesting (though probably not easy) to study the analogous problem in the random setting, estimating the average efficiency. This subject is not completely new (compare Section 6 of [LM]).
\smallskip

The number of open questions can be increased further if we look at different sorts of analogues of the Steiner symmetrization such as Schwarz, Blaschke and Minkowski (see [SY] for the definitions) and the type of space, considering, as suggested by Mani-Levitska, besides the Euclidean, also spherical and hyperbolical spaces.
\bigskip

\noindent {\bf Acknowledgements}
\bigskip

The author wants to express his thanks to Ingrid Carbone for her critical reading of the manuscript, which improved the presentation.
\bigskip\bigskip

\bigskip

Universit\`a della Calabria

Dipartimento di Matematica

87036 Arcavacata di Rende (CS) - Italia

Ponte Bucci, cubo 30B

volcic@unical.it
\bigskip\bigskip

\noindent {\bf References}
\bigskip

\noindent [AFP] F. Ambrosio, N. Fusco, D. Pallara, {\it Functions of bounded variations and free discontinuity problems}. Oxford University Press, Oxford, 2000.

\smallskip
\noindent [BG] G. Bianchi, P. Gronchi, Steiner symmetrals and their distance from a ball,
{\it Israel J. Math.} {\bf 135} (2003), 181-192. 

\smallskip
\noindent [BLM] J. Bourgain, J.  Lindenstrauss, V. Milman, {\it Estimates related to Steiner symmetrizations, in Geometric Aspects of Functional Analysis}, eds J. Lindenstrauss and V. Milman, Springer, Lecture Notes in Math. {\bf1376} (1989), 264-273.

\smallskip

\noindent [C] R. Caccioppoli, Elementi di una teoria generale dell'integrazione {\it k}-dimen\-sionale in uno spazio {\it n}-dimensionale, {\it Atti del Quarto Congresso dell'Unione Matematica Italiana, Taormina, 1951, vol. II},  pp. 41--49. Casa Editrice Perrella, Roma, 1953.
\smallskip

\noindent [Ch] K. L. Chung, {\it A Course in Probability Theory. Third edition}. Academic Press, Inc., San Diego, CA, 2001.
\smallskip

\noindent [DG ] E. De Giorgi, Su una teoria generale della misura $(r-1)$-dimensiona\-le in uno spazio ad $r$ dimensioni, {\it Ann. Mat. Pura Appl. (4)} {\bf 36}  (1954), 191-213.

\smallskip
\noindent [Ga] R. J. Gardner,  {\it Geometric tomography. Second Edition.} Encyclopedia of Mathematics and its Applications {\bf 58}, Cambridge University Press, Cambridge, 2006.
\smallskip

\noindent [G] W. Gross, Die Minimaleigenschaft der Kugel, {\it Monatsh. Math. Phys.}  {\bf 28}   (1917), no. 1,  77-97.
\smallskip

\noindent [Gr] P. M. Gruber, The space of convex bodies, in {\it Handbook of Convex Geometry}, ed. by P. M. Gruber and J. M. Wills, Horth-Holland, Amsterdam (1993), 301-318.
\smallskip

\noindent [H] H. Hadwiger, Einfache Herleitung der isoperimetrischen Ungleichung f\"ur abgesclossene Punktmengen, {\it Math. Ann.} {\bf 124} (1952), 158-160.
\smallskip

\noindent [KM1] B. Klartag, V. Milman, Isomorphic Steiner symmetrizations, {\it Invent. Math.}  {\bf 153} (2003),  no. 3, 463-485.
\smallskip

\noindent [KM2] B. Klartag, V. Milman, Rapid Steiner symmetrization of most of a convex body and the slicing problem, {\it Combin. Probab. Comput.}  {\bf 14} (2005), no. 5-6, 829-843. 
\smallskip

\noindent [KN] L. Kuipers, H. Niderreiter,  {\it Uniform distribution of
sequences. Pure and Applied Mathematics}. Wiley-Interscience, New York-London-Sidney, 1974.

\smallskip

\noindent [LM] J. Lindenstrauss, V. Milman, The local theory of normed spaces and its applications to convexity,  in {\it Handbook of Convex Geometry}, ed. by P. M. Gruber and J. M. Wills, Horth-Holland, Amsterdam (1993), 1149-1220.
\smallskip

\noindent [M] P. Mani-Levitska, Random Steiner symmetrizations, {\it Studia Sci. Math. Hungar.} 
{\bf 21} (1986), no. 3-4, 373-378. 
\smallskip

\noindent [S] J. Steiner, Einfache Beweis der isoperimetrischen Haupts\"atze, {\it J. reine angew. Math.} {\bf 18} (1838), 281-296.
\smallskip

\noindent [SY] J. R. Sangwine-Yager, Mixed volumes, in {\it Handbook of Convex Geometry}, ed. by P. M. Gruber and J. M. Wills, Horth-Holland, Amsterdam (1993), 43-71.
\smallskip

\noindent [T] G. Talenti, The standard isoperimetric theorem, in {\it Handbook of Convex Geometry}, ed. by P. M. Gruber and J. M. Wills, Horth-Holland, Amsterdam (1993), 73-123.
\smallskip

\end{document}